\documentclass{amsart}
\usepackage{mathrsfs}
\usepackage{hyperref}

\usepackage{color}

\usepackage{marginnote}

\def\R{{\mathbb R}}

\newtheorem{theorem}{Theorem}[section]
\newtheorem{proposition}[theorem]{Proposition}
\newtheorem{corollary}[theorem]{Corollary}
\newtheorem{definition}[theorem]{Definition}

\def\End{\mathrm {End}}

\def\ad{\mathrm{ad}}

\def\di{\displaystyle}

\begin{document}
\title[The Ricci flow  for circle bundles over surfaces]{The Ricci flow  for circle bundles over surfaces}

\author{Arash Bazdar}
\address{Department of Mathematics, Xiamen University Malaysia, 43900 Selangor, Malaysia}
\email{bazdar.arash@xmu.edu.my}

\author{Georgios Fotopoulos}
\address{Academic Support Department, Abu Dhabi Polytechnic, P.O. Box 111499, Abu Dhabi, United Arab Emirates}
\email{georgios.fotopoulos@actvet.gov.ae}

\begin{abstract}
In this work, we study and solve the normalized Ricci flow equation for circle bundles over surfaces. Moreover, we  study  the asymptotic behavior of the solutions and their connections to some model geometries.
\end{abstract}


\maketitle
	
\section{Introduction}
The Ricci flow equation, introduced by Hamilton in his seminal paper \cite{H82}, has proven to be a powerful tool in the study of geometric structures on manifolds. Hamilton's objective was to utilize the Ricci flow to prove Thurston’s geometrization conjecture, which states that every closed manifold admits a geometric decomposition.

In this paper, we study the Ricci flow for circle bundles over surfaces, endowed with  locally homogeneous $\mathbb S^1$-connections. These connections arise in the study of quantum field theory and gauge theory (see \cite{bleecker}), and a thorough  understanding of their geometric properties is crucial for a deeper comprehension of these physical theories. 
 It is worth noting that the concept of a circle bundle can be extended to the notion of Seifert fibration, which originally appeared in the classification of spherical $3$-manifolds (see \cite{Sc,Seifert33}).
Moreover, it has been shown that a compact orientable $3$-manifold is a Seifert 3-manifold if and only if it admits a geometry modeled on one of the following six geometries: $\mathbf{S}^3$, $\mathbf{S}^2\times \R$, $\mathbf{H}^2\times\R$, $\widetilde{\mathbf{SL}_2(\R)}$, $\mathbf{Nil}$ and $\mathbf{Sol}$ (see \cite[Chap. 1 and 2]{Sc,mailot}). In particular, we show that the solutions of the normalized Ricci flow for circle bundles over surfaces give rise to the aforementioned six geometries. 

The Ricci flow has emerged as a crucial tool in resolving both the Poincaré and Thurston conjectures, highlighting its significance in the study of geometric manifolds (see \cite{an,mf,s2}). 
Research on the Ricci flow applied to specific bundles over surfaces has yielded significant results.  In a recent preprint \cite{Yeroshkin21}, the author studied the Ricci flow on torus bundles, while a similar problem but in a different setting was studied in \cite{S10,s2}, where the author explored the Ricci Yang-Mills flow, aiming to simplify its behavior compared to the Ricci flow. This investigation  introduced connections with specific curvature properties, with the expectation that these modifications would lead to a flow exhibiting simpler characteristics than the Ricci flow.  In \cite{S22} the author gave a comprehensive  account of the global existence and convergence of the Ricci Yang-Mills flow on torus bundles over Riemann surfaces. For a more detailed exploration of Ricci Yang-Mills flow, see \cite{y1,y2,y3}.

 In our approach, we fix a Yang-Mills connection on the circle bundle in order to simplify the system, while allowing the size of the fiber and the metric on the base manifold to evolve over time. In contrast, the Ricci Yang-Mills flow allows both the connection and the metric on the base to evolve over time under the additional hypothesis that the size of the fiber remains fixed \cite{y3}.
Motivated by the works in \cite{Ba, IJ}, where the former explores geometric principal bundles arising from locally homogeneous connections on bundles over geometric base manifolds, and the latter investigates the Ricci flow of all locally homogeneous Riemannian metrics on closed 3-manifolds by analyzing the Ricci flow of their corresponding homogeneous models, we have undertaken a study of the Ricci flow on circle bundles over surfaces. In particular, we would like to highlight the similarities between our results and those in \cite{IJ}, where their approach involves examining both Bianchi and non-Bianchi classes of geometries, while our approach focuses on metrics obtained on principal bundles by considering a connection, along with the metrics on the base manifold and the structure group of the bundle. If the connection and the metric on base are locally homogeneous, the induced metric on the bundle is also locally homogeneous, thereby defining a geometric structure on the bundle. This wide range of geometric manifolds makes the study of their Ricci flow an engaging and valuable research problem. An especially promising direction for future research involves examining the $\mathrm{SU}(2)$-bundle over a surface endowed with a locally homogeneous connection, leading to the construction of five-dimensional geometric manifolds (see \cite{Ba2} for further details).

The paper is structured as follows. In Section \ref{sec2} we introduce the necessary background on Levi-Civita connection, Riemann and Ricci curvature tensor for a circle bundle over surfaces endowed with a connection. These foundational calculations are essential for deriving the Ricci flow equation and are also significant in their own right, as they are not readily available in the  existing literature. In Section \ref{sec3}, we present various results concerning the existence of solutions to the Ricci flow equation and their geometric properties. We demonstrate that the normalized Ricci flow equation for a circle bundle reduces to a system of nonlinear, first-order ordinary differential equations (ODEs). Under certain conditions, this system can be simplified and solved explicitly. Specifically, for the curvature two-form of the connection, we investigate different scenarios, each corresponding to a system of ODEs, which we study in detail. We show that the circle bundle can exhibit one of the six geometries from Thurston’s list, and we provide the solutions along with their asymptotic behavior. We conclude with Section \ref{sec4}, which summarizes our work. 

\section{Locally homogeneous \texorpdfstring{$\mathbb S^1$}{S\^1}-connections over  surfaces}\label{sec2}

\subsection{Metrics on principal \texorpdfstring{$K$}{K}-bundles} \label{TheProblem} 
A Riemannian metric $g$ on a differentiable manifold $M$ is said to be locally homogeneous if, for any two points $x$ and $y$ in $M$, there exist open neighborhoods $U$ containing $x$ and $V$ containing $y$ and an isometry $\varphi:U\to V$ such that $\varphi(x)=y$. Motivated by this definition, we define a locally homogeneous connection as follows. 

\begin{definition}\label{LHTDef}
Let $(M, g)$ be a connected, locally homogeneous Riemannian manifold, let $K$ be a connected, compact Lie group, and let $\pi : P \to M$ be a principal $K$-bundle on $M$. A connection $A$ on $P$ is said to be locally homogeneous if, for any two points $x$ and $y$ in $M$, there exists an isometry $\varphi:U\to V$ between open neighborhoods $U$ containing $x$ and $V$ containing $y$ with $\varphi(x)=y$, and an $\varphi$-covering bundle isomorphism $\Phi:P_U \to P_V$ such that $\Phi^*(A_V) = A_U$.
\end{definition}

Note that, the subscripts $U$ and $V$ denote the restriction of the indicated objects to the open subsets $U, V \subset M$. Following the notation used in \cite{KN}, we recall that $A$ stands for a $K$-invariant horizontal distribution of $P$. Let $g_K$ be an $\ad$-invariant inner product on $\mathrm{Lie}(K)$. We endow the total space $P$ with a Riemannian metric $g_A$, called the bundle metric, by fixing an $\ad$-invariant inner product on the Lie algebra of $K$, characterized by the following conditions:

\begin{enumerate}
\item The canonical bundle isomorphism $V_P\simeq P\times \mathrm{Lie}(K) $ is an orthogonal bundle isomorphism with respect to the inner products defined by $g_A$ and $g_K$.
\item The restriction of the differential map $\pi_* : T_P\to T_M$ to the horizontal sub-bundle $A\subset T_P$ gives an orthogonal bundle isomorphism $A\to \pi^*(T_M)$.
\item The direct sum decomposition of the tangent bundle $T_P\simeq A\oplus V_P$ is $g_A$-orthogonal. 
\end{enumerate}
If the connection $A$ is locally homogeneous, the bundle metric $g_A$ on $P$ is locally homogeneous as well.

Next, we recall a well-known theorem of Singer \cite{Si}, which asserts that any complete, simply connected, locally homogeneous Riemannian manifold is homogeneous. Let $p:\widetilde M\to M$ be the universal cover of $M$, then the metric $\tilde g:=p^*g$ defines a complete locally homogeneous Riemannian metric on the simply connected manifold $\widetilde M$. Hence, by Singer's theorem, $\tilde{g}$ is homogeneous. In particular, $\widetilde M$ admits a connected transitive group of isometries $G$. Therefore, for a locally homogeneous Riemannian manifold $(M,g)$, there exists a homogeneous model geometry given by $(\widetilde M,G)$, where a model geometry $(X,G)$ is a simply connected manifold $X$ together with a transitive action of a Lie group $G$ on $X$ with compact stabilizers. Generalizing the method introduced by Singer, one can prove a similar classification theorem for locally homogeneous connections (see \cite{Ba}).

Recall that for three-dimensional manifolds, Thurston \cite{Th} classified eight maximal model geometries with compact quotients, namely $\mathbf{E}^3$, $\mathbf{S}^3$, $\mathbf{H}^3$, $\mathbf{S}^2\times \R$, $\mathbf{H}^2\times\R$, $\widetilde{\mathbf{SL}_2(\R)}$, $\mathbf{Nil}$ and $\mathbf{Sol}$ (see also \cite{Sc}). A geometric structure on a manifold $M$ is a diffeomorphism from $M$ to $X/\Gamma$ for some model geometry $(X,G)$, where $\Gamma\subset G$ is a discrete subgroup of $G$ acting freely on $X$. For any of the Thurston geometries, except $\mathbf{H}^3$ and $\mathbf{Sol}$, there exists a compact $3$-manifold with a geometric structure that is associated with locally homogeneous $\mathbb S^1$-connections. For instance, any non-trivial $\mathbb S^1$-bundle on a surface of genus $g$ has a geometric $(X,G)$-structure, where 

\begin{enumerate}
\item $(X,G)=\mathbf{S}^3$, when $g=0$,
\item $(X,G)=\mathbf{Nil}$, when $g=1$,
\item $(X,G)=\widetilde{\mathbf{SL}_2(\R)}$, when $g\geq 2$.
\end{enumerate}

The trivial $\mathbb S^1$-bundles over surfaces have geometric structures with model geometry $\mathbf{S}^2\times \R$, $\mathbf{E}^3$, or $\mathbf{H}^2\times\R$. Moreover, any $\mathbb S^1$-bundle over a surface has a metric, which is associated with a locally homogeneous connection in the sense of Definition \ref{LHTDef}. These examples show the abundance of geometric principal $\mathbb S^1$-bundles in the class of geometric $3$-manifolds. 
 
\subsection{The Levi-Civita connection associated with a bundle metric on circle bundles over surfaces } \label{sec12}Let $\pi:P\to M$ be a principal $\mathbb S^1$-bundle on a surface $M$, let $g$ be a Riemannian metric on $M$, and let $A$ be a connection on $P$. We identify the Lie algebra of $\mathbb S^1$ with $\R$ and denote the connection $1$-form by $\omega_A\in A^1(P,\R)$ and endow the total space $P$ with the bundle metric $g_A$. Let $\varphi=(\varphi _1,\varphi _2)$ be an orthonormal frame of a surface $(M,g)$ over an open neighborhood $U\subset M$ and let $(e_1:=\tilde \varphi_1,e_2:=\tilde \varphi_2)$ be the $A$-horizontal lifts of $\varphi$ on $P_U:=\pi^{-1}(U)$. Therefore, we have
\[ \omega(e_i)=0 \ \text{,} \ \pi_*(e_i)=\varphi_i \ \quad \text{for} \quad i=1,2. \]
Let $\xi$ be an element of the Lie algebra of $\mathbb S^1$. The vector field $\xi^\#$ defined by
\[ \xi^\#_y:=\frac{d}{dt}\Big|_{t=0} \big(y\exp(t\xi)\big), \]
is called the fundamental vector field corresponding to $\xi$. For any point $y$ in the total space $P$, we define 
\[ f(y):=\sqrt{g_A(\xi^\#_y,\xi^\#_y)}, \] 
which is a constant since the metric on the Lie algebra of $\mathbb S^1$ is fixed. We define the section $e_3:P\to V_P$ by
\[ e_3(y) = \frac{1}{f(y)}\xi_y^\# ,\quad \text{for } y\in P. \]
Evidently, $(e_1,e_2,e_3)$ is an orthonormal frame on $P_U$ with respect to the metric $g_A$. Let $(\eta^1, \eta^2)$ be the dual frame associated with $(\varphi _1,\varphi _2)$ over $U\subset M$. Setting $\theta^1=\pi^*\eta^1$, $\theta^2=\pi^*\eta^2$ and $\theta^3=f\omega$, we obtain a co-frame $(\theta^1,\theta^2,\theta^3)$ dual to $(e_1,e_2,e_3)$ on $P_U$.

By a well-known result of differential geometry, see \cite[Proposition 5.32]{Mo}, there exists a unique skew-symmetric matrix $(\eta^i_j)$ of real-valued $1$-forms $\eta^i_j\in A^1(U,\R)$ such that 
\begin{equation} \label{sM1}
d\eta ^1 = -\eta ^1_2\wedge \eta ^2 \quad\text{and}\quad d\eta ^2 = -\eta ^2_1\wedge \eta ^1.
\end{equation}
Taking the pullback on both sides of the equations \eqref{sM1}, we obtain
\begin{equation} \label{sM2}
d\theta^1 = -\pi^*\eta ^1_2\wedge \theta ^2 \quad\text{and}\quad d\theta^2 = -\pi^*\eta ^2_1\wedge \theta ^1.
\end{equation}
The curvature $2$-form of the connection $A$ is denoted by $F_A =d\omega$. Using the dual frame $(\theta^1,\theta^2,\theta^3)$ on $P_U$, the curvature $2$-form $F_A$ can be written as
\begin{equation*}
F_A=c\,\theta^1\wedge\theta^2 , \quad c\in \mathcal C^\infty(P_U,\R).
\end{equation*}
Differentiating $\theta^3$, since $f$ is constant, we get
\begin{equation}\label{dteta3}
d\theta^3=d(f\omega)=fd\omega+df\wedge \omega= fc\,\theta^1\wedge \theta^2.
\end{equation}
The equations (\ref{sM2}) and (\ref{dteta3}) lead to the following system of differential equations
\begin{align} \label{sP0}
d\theta^1 &= -\pi^*\eta ^1_2\wedge \theta ^2,\nonumber\\ 
d\theta^2 &= -\pi^*\eta ^2_1\wedge \theta ^1,\\
d\theta^3 &= fc\, \theta^1\wedge \theta^2.\nonumber
\end{align}
Our aim is to find a skew-symmetric matrix $(\theta^i_j)$ of real-valued $1$-forms on $P$ which satisfy the differential equations
\begin{equation}\label{sp1}
d\theta^i=-\sum_j\theta^i_j\wedge\theta^j ,
\end{equation} 
where the solutions of these equations are the components of the Levi-Civita connection $1$-form associated with the metric $g_A$. In the following proposition we compute the aforementioned solutions.

\begin{proposition} \label{LC}
Let $\pi:P\to M$ be a principal $\mathbb S^1$-bundle over a surface $(M,g)$ endowed with a connection $A$. Using the same notation as above, the Levi-Civita connection $(\theta^i_j)$ of the bundle metric $g_A$ is given by
\begin{align}
\theta^1_2 & =\pi^*(\eta^1_2)-\frac{1}{2}f c\, \theta^3, \\
\theta^1_3 & =-\frac{1}{2}f c\, \theta^2, \\
\theta^2_3 & =\frac{1}{2}f c\, \theta^1.
\end{align}
\end{proposition}

\begin{proof}
Using the system \eqref{sP0}, we find the real-valued $1$-forms $\mu^i_j$ on $P_U$ which satisfy the system \eqref{sp1}, that is,
\[ 
\mu^i_j:=\begin{pmatrix}
0&\pi^*\eta^1_2 &0 \ \\ 
\pi^*\eta^2_1 & 0 &0 \ \\ 
\ fc\,\theta^2& 0 & 0\ 
\end{pmatrix}.
\]
Since this solution is not skew-symmetric, we modify it in an appropriate way as follows. We set
\begin{equation}\label{solotionLC}
\theta^i_j:=\mu^i_j+{\sum}_k S^i_{jk}\theta^k,
\end{equation}
where $S$ is a unique tensor (to be calculated) such that $S^i_{jk}=S^i_{kj}$ for all $i,j,k$. Then, it is easy to see that $\theta^i_j$ satisfy \eqref{sp1}. For the modified solution
\[ \theta^i_j:=\mu^i_j+ {\sum}_k S^i_{jk}\theta^k, \] 
to be skew-symmetric, the following conditions need to be satisfied: 
\begin{align}
\theta^i_j&=-\theta^j_i,\label{MSys2}\\
 S^i_{jk} &=S^i_{kj}\label{MSys3}.
\end{align}
Note that, for $i=j$, \eqref{MSys2} yields
\begin{equation}
0=\theta^i_i=\mu^i_i+{\sum}_kS^i_{ik}\theta^k={\sum}_k S^i_{ik}\theta^k.
\end{equation}
So, we have that 
\[ S^i_{ik}=S^i_{ki}=0, \quad\text{for all } i,k, \]
where we used \eqref{MSys3}. Next, for $i \ne j$ and taking all possible combinations, we derive 
\begin{align*}
&S^1_{2k} = -S^2_{1k}, \quad\text{for all } k,\\
&S^3_{11} =S^1_{33}=0,\\
&S^3_{12} + S^1_{32}=-fc,\\
&S^3_{22} =S^2_{33}=0,\\
&S^3_{12} = S^1_{32} =-\frac{1}{2}fc,\\ 
&S^2_{13} =\frac{1}{2}fc.
\end{align*}
Substituting $S^i_{jk}$ and $\mu_i^j$ back into \eqref{solotionLC}, we obtain the desired result.
\end{proof}

\subsection{The curvature tensor of a bundle metric on \texorpdfstring{$\mathbb S^1$}{S\^1}-bundles over  surfaces}
Let $\pi:P\to M$ be a principal $\mathbb S^1$-bundle over a  surface $(M,g)$ endowed with a bundle metric $g_A$ and let $(\theta^i_j)$ be the Levi-Civita connection of $g_A$ obtained in Proposition \ref{LC}. The curvature $2$-form is denoted by $\Omega^\theta$ and its components $(\Omega^i_j)$ are given in the following proposition.

\begin{proposition}\label{CurvP} Let $\theta=(\theta^i_j)$ be the connection $1$-form of the Levi-Civita connection of $g_A$ and $K$ be the Gauss curvature of the surface $(M,g)$. Then, the components of the curvature $2$-form $\Omega^\theta=(\Omega^i_j)$ are given by
\begin{align*}
\Omega_{2}^{1} & = \Big( \pi^*K-\frac{3}{4}f^2c^2 \Big) \theta^{1}\wedge \theta^{2}-\frac{1}{2}f\alpha\theta^{1}\wedge \theta^{3}-\frac{1}{2}f\beta\theta^{2}\wedge \theta^{3},\\
\Omega_{3}^{1} &=-\frac{1}{2}f\alpha\theta^1\wedge \theta^2+\frac{1}{4}f^2c^2\theta^1\wedge\theta^3,\\
 \displaystyle \Omega_{3}^{2}& =-\frac{1}{2}f\beta \theta^1\wedge\theta^2+\frac{1}{4}f^2c^2\theta^2\wedge\theta^3,
\end{align*}
where $d c=\alpha\theta^1+\beta\theta^2$ and $\alpha, \beta$ are smooth functions on $P_U$.
\end{proposition}

\begin{proof}
Using the structure equation 
\[ \Omega^i_j=d\theta^i_j+\sum_k\theta^i_k\wedge\theta^k_j , \]
we determine the components of $\Omega^\theta$. Starting with the computation of $\Omega_{2}^{1}$, we have 
\[ \Omega_{2}^{1}=d\theta_{2}^{1}+\theta_{3}^{1}\wedge \theta_{2}^{3}. \]
Using Proposition \ref{LC}, the first term of the right-hand side is given by
\begin{align}\label{Cteta21}
d\theta_{2}^{1} & =\pi^{*}d\eta_{2}^{1}-\frac{1}{2} d(fc\, \theta^{3})\nonumber\\ 
&= \pi^*K(\theta^1\wedge\theta^2)-\frac{1}{2} \big(f( dc\wedge \theta^{3})+c\,(df\wedge\theta^{3} )+fc\,d\theta^{3}\big)\nonumber\\
&=\pi^*K(\theta^1\wedge\theta^2) -\frac{1}{2}\big(f\alpha\, \theta^{1}\wedge\theta^{3}+f\beta\theta^{2}\wedge\theta^{3}+f^2c^2\theta^{1}\wedge \theta^{2}\big),
\end{align}
where $K$ denotes the Gauss curvature of the surface $(M,g)$ and we used 
\begin{align*}
d c &=\alpha\theta^1+\beta\theta^2,\\
d\theta^{3} &= d(f\omega)=fd\omega+df\wedge\omega =fc\, \theta^{1}\wedge \theta^{2}.
\end{align*}
The second term is given by
\begin{equation}\label{delta1323}
\theta_{3}^{1}\wedge \theta_{2}^{3} = \Bigl( -\frac{1}{2}fc\theta^{2} \Bigr) \wedge \Bigl( -\frac{1}{2}fc\theta^{1} \Bigr) = -\frac{1}{4}f^{2}c^{2}\theta^{1}\wedge \theta^{2}.
\end{equation}
Therefore, by (\ref{Cteta21}) and (\ref{delta1323}), we get
\[ \Omega_{2}^{1}= \Big(\pi^* K-\frac{3}{4}f^2c^2\Big)\theta^{1}\wedge \theta^{2}-\frac{1}{2}f\alpha\theta^{1}\wedge \theta^{3}-\frac{1}{2}f\beta\theta^{2}\wedge \theta^{3}. \]
Moving on to the computation of $\Omega_{3}^{1}$, we have
\[ \Omega_{3}^{1}=d\theta_{3}^{1}+\theta_{3}^{2}\wedge \theta_{2}^{1}. \]
Using Proposition \ref{LC}, we have 
\begin{align}\label{dfcteta2}
d\theta^1_3 &= d \Bigl(-\frac{1}{2}fc\,\theta^2\Bigr) =-\frac{1}{2} ( fdc\wedge\theta^2+fc\,d\theta^2 ) \nonumber\\
& =-\frac{1}{2} \big( f\alpha\, \theta^{1}\wedge \theta^{2}+fc (-\pi ^*\eta^2_1\wedge\theta^1)\big).
\end{align}
We also have
\begin{align}\label{teta2132l}
\theta_{2}^{1}\wedge \theta_{3}^{2} &=\Big(\pi^{*}\eta_{2}^{1}-\frac{1}{2}fc\,\theta^{3} \Big)\wedge \Big(\frac{1}{2}fc\,\theta^{1} \Big) \nonumber\\
&=\frac{1}{2}fc\, (\pi^*\eta^1_2)\wedge\theta^1+\frac{1}{4}f^2c^2\theta^1\wedge\theta^3.
\end{align}
Adding \eqref{dfcteta2} and \eqref{teta2132l}, we get
\[ \Omega_{3}^{1}=-\Big(\frac{1}{2}f\alpha\Big) \theta^1\wedge \theta^2+\Big(\frac{1}{4}f^2c^2\Big)\theta^1\wedge\theta^3 . \]
Finally, we determine
\[ \Omega_{3}^{2}=d\theta^2_3+\theta_{1}^{2}\wedge \theta_{3}^{1}, \]
where the term $d\theta_{3}^{2}$ is given by 
\begin{align}\label{dteta32l}
d\theta_{3}^{2} & =\frac{1}{2}d(fc\theta^1)=\frac{1}{2}\left( fc\,d\theta^1+fdc\wedge\theta^1 \right)\nonumber \\
&=-\frac{1}{2} fc\,(\pi^*\eta^1_2)\wedge\theta^2 - \frac{1}{2}f\beta\theta^1\wedge\theta^2,
\end{align}
and the term $\theta_{1}^{2}\wedge \theta_{3}^{1}$ is given by
\begin{equation}\label{teta1231l}
\theta_{1}^{2}\wedge \theta_{3}^{1}=(-\pi^{*}\eta_{2}^{1}+\frac{1}{2}fc\,\theta^{3} )\wedge (-\frac{1}{2}fc\,\theta^{2} )=\frac{1}{2}fc\,(\pi^*\eta^1_2)\wedge\theta^2+\frac{1}{4}f^2c^2\theta^2\wedge\theta^3. 
\end{equation}
Combining \eqref{dteta32l} and \eqref{teta1231l} we obtain
\[ \Omega_{3}^{2}=-\frac{1}{2}f\beta\theta^1\wedge\theta^2+\frac{1}{4}f^2c^2\theta^2\wedge\theta^3. \qedhere \]
\end{proof}

\subsection{The Ricci curvature of a metric bundle on principal \texorpdfstring{$\mathbb S^1$}{S\^1}-bundles over  surfaces}
Let $\pi: P\to M$ be a principal $K$-bundle over $M$, endowed with the bundle metric $g_A$. Let $(e_i)$ be a $g_A$-orthonormal frame on $P_U$, where $U$ is an open subset of $M$. For an arbitrary point $y\in P_U$ and two vectors $v,w\in T_yP$, the Ricci curvature tensor is given by 
\[ \mathrm{Ric}(v,w)=\sum_m g_A\left(R^\theta(v,e_m)e_m,w\right), \]
where $R^\theta\in A^2(P,\End{(T_P)})$ is the Riemann curvature tensor associated to the Levi-Civita connection $\theta$. Hence, the components of the Ricci curvature tensor in the $g_A$-orthonormal frame $(e_i)$ are given by
 \[ \mathrm{R}_{lk}=\mathrm{Ric}(e_l,e_k)=\sum _mg_A\left( R^\theta(e_l,e_m)e_m,e_k\right). \]
The Riemann curvature tensor $R^\theta$ is related to the curvature $2$-form $\Omega^\theta=(\Omega^i_j)$ by 
\[ R^\theta(e_l,e_m)e_m=\sum_{i}\Omega^i_m(e_l,e_m)e_i. \] 
So, we obtain
\begin{equation}\label{Riclk}
\mathrm{R}_{lk}=\sum_{m,i}g_A\left( \Omega^i_m(e_l,e_m)e_i,e_k \right)=\sum_m \Omega^k_m(e_l,e_m).
\end{equation}
In the proposition below, we compute the components of the Ricci tensor of $g_A$, using the Proposition \ref{CurvP} and the equations \eqref{Riclk}. 
\begin{proposition}\label{RICCh}
	Let $\pi: P\to M$ be a principal $\mathbb S^1$-bundle on $M$ endowed with a connection $A$, let $\omega\in A^1(P,\mathbb{R})$ be the connection 1-form, and let $F_A\in A^2(M)$ be the curvature form of $A$. The Ricci curvature tensor of $g_A$ is given by
 \[ \mathrm{Ric}(g_A)= \Big(\pi^*K-\frac{1}{2}f^2\pi^*|F_A|^2_g \Big)\,\pi^*g-\frac{1}{2}f^2\pi^*(d^*F_A)\otimes\omega +\frac{1}{2}f^4\pi^*|F_A|^2_g\, \omega \otimes \omega, \]
where $K$ is the Gauss curvature of the base manifold and $d^*$ is the co-differential operator.
\end{proposition}

\begin{proof}
The components $\mathrm{R}_{1k}, k=1,2,3$ are given by
\[ \mathrm{R}_{1k}=\Omega^k_2(e_1,e_2)+\Omega^k_3(e_1,e_3), \]
which gives
\[ \mathrm{R}_{11} = \pi^*K-\frac{1}{2}f^2c^2, \quad \mathrm{R}_{12}=0, \quad \mathrm{R}_{13}=\frac{1}{2}f\beta. \]
In a similar way, we obtain
\[ \mathrm{R}_{21}=0, \quad \mathrm{R}_{22}=\pi^*K-\frac{1}{2}f^2c^2, \quad \mathrm{R}_{23}=-\frac{1}{2}f\alpha, \]
and
\[ \mathrm{R}_{31}=\frac{1}{2}f\beta, \quad \mathrm{R}_{32}=-\frac{1}{2}f\alpha, \quad \mathrm{R}_{33}=\frac{1}{2}f^2c^2. \]
Therefore, the Ricci tensor of the metric $g_A$ with respect to the frame $(e_1,e_2,e_3)$ can be represented by
\[ \mathrm{Ric}(g_A)=\frac{1}{2} \begin{pmatrix}
2\pi^*K-f^2c^2&0 & f\beta \\ 
0&2\pi^*K-f^2c^2 &-f\alpha\\
-f\alpha& f\beta&f^2c^2
\end{pmatrix}. \]
Recall that if $\omega \in A^1(P,\mathbb{R})$ denotes the connection 1-form, then the curvature form satisfies
\[ \pi^*F_A = d\omega = c \,\theta^1 \wedge \theta^2, \]
and it is easy to see that $\pi^*|F_A|^2_g=c^2$. Finally, taking into account the identity 
\[ -\frac{1}{2}f^2\pi^*(d^*F_A)\otimes\omega =\frac{1}{2}(f\beta)\theta^1\otimes\theta^3-\frac{1}{2}(f\alpha)\theta^2\otimes \theta^3, \]
we obtain the desired result.
\end{proof}

For the special case where $A$ is a Yang-Mills connection, that is, $d^*F_A=0$, we have the following result.

\begin{corollary} \label{scalarwithfconst}The following hold.
\begin{enumerate}
\item If $R_g=2K_g$ denotes the scalar curvature of the metric $g$, then the scalar curvature of $g_A$ is given by 
\[ R_{g_A}=\pi^*R_g-\frac{1}{2}f^2\pi^*|F_A|^2_g. \]
\item The horizontal and vertical subspaces are orthogonal relative to the Ricci tensor of the metric $g_A$ and we have
\[ \mathrm{Ric}(g_A) = \Big( \pi^*K-\frac{1}{2}f^2\pi^*|F_A|^2_g \Big) \pi^*g+\frac{1}{2}f^4\pi^* |F_A|^2_g\, \omega \otimes \omega . \]
\end{enumerate}
\end{corollary}

\section{The Ricci flow equation of a metric bundle on principal \texorpdfstring{$\mathbb S^1$}{S\^1}-bundles over  surfaces}\label{sec3}
\subsection{The Ricci flow equation on principal \texorpdfstring{$\mathbb S^1$}{S\^1}-bundles over  surfaces}\label{sec21}
Let $(M,g)$ be a connected, oriented  surface and let $P$ be a $\mathbb S^1$-bundle over $M$, endowed with a connection $A$. The bundle metric on $P$ is defined by 
\begin{equation}\label{metricgA}
	g_A=\pi^*g+f^2\omega\otimes\omega,
\end{equation} 
where $\omega$ is the connection 1-form and $f$ is a constant as in Section \ref{sec12}. 

Since we are interested in studying the Ricci flow equation of the metric $g_A$, we consider \eqref{metricgA} with respect to time $t$. Therefore, differentiating with respect to time, we get
\[ \dot{g}_A(t)=\pi^*\dot{g_t}+f_t^2\dot{\omega_t}\otimes\omega_t+f_t^2\omega_t\otimes\dot{\omega_t}+2f_t\dot{f_t}\,\omega_t\otimes \omega_t . \]
Recall that the Ricci flow equation of the metric $g_A$ is given by
$$\dot{g}_A(t)=-2\mathrm{Ric}(g_A(t)).$$
Therefore, by the above equations and Proposition \ref{RICCh}, the Ricci flow equation of $g_A$ is equivalent to the following system of first-order differential equations:
\begin{align}\label{RFY}
\dot{g_t}&= -2\mathrm{Ric}(g_t)+f_t^2\left |F_{A_t} \right |^2_{g_t}g_t,\\ 
\dot{\omega_t }&= d^*F_{A_t},\\ 
\dot{f_t}&= -\frac{1}{2}f_t^3\left |F_{A_t} \right |^2_{g_t}.
\end{align}

Next, we restrict ourselves to the case where the base manifold has constant curvature. To this end, we assume $(M,g_0)$ is a connected  surface with constant Gauss curvature $K_0$ and $A_0$ be a Yang-Mills connection on $P$. 

Setting $g_t=\lambda_tg_0$ and $\omega_t=\omega_0$, where $\omega_0$ is the connection 1-form of $A_0$, the system \eqref{RFY} becomes 
\begin{align} \label{rfy1}
		\dot{\lambda_t}& = -2K_0+F_0 f_t^2\lambda_t^{-1},\\
		\dot{f_t}& = -\frac{1}{2}F_0 f_t^3\lambda_t^{-2},
\end{align}
where we used the fact that the Ricci tensor is scale invariant, that is,
\begin{equation*}\label{manbaror231}
	\mathrm{Ric}(g_t)=\mathrm{Ric}(\lambda_tg_0)=\mathrm{Ric}(g_0)=K_0g_0
\end{equation*}
and 
\begin{equation*}\label{manbaror}
\left |F_{A_0} \right |^2_{g_t}=\lambda_t^{-2}\left |F_{A_0} \right |^2_{g_0}=\lambda_t^{-2}F_0,\quad F_0:=\left |F_{A_0} \right |^2_{g_0}.
\end{equation*}
Note that $|\cdot|_g$ denotes the norm induced by the metric $g$ on 2-forms.

\subsection{Solutions of the normalized Ricci flow equation for locally homogeneous \texorpdfstring{$\mathbb S^1$}{S\^1}-connections on  surfaces} 
Let $(M,g)$ be a connected, oriented,  surface endowed with a metric with a constant curvature, let $P$ be a principal $\mathbb S^1$-bundle on $M$, and let $A$ be a connection on $P$. Note that $A$ is a locally homogeneous connection if and only if $A$ is Yang-Mills (see \cite{Ba}).

Next, we consider the normalized Ricci flow equation of the metric $g_A$ on $P$, that is,
\begin{equation}\label{NRFEq}
\frac{\partial}{\partial t}g_A(t)=-2\mathrm {Ric} \ g_A(t)+\frac{2}{3}r(t)g_A(t),
\end{equation}
where 
\[ r(t):= \frac{\di\int_PR(t)d \mu_t}{\di\int_P d \mu_t} \]
is the average of the scalar curvature $R(t)$ of $g_A(t)$ and where $\mu_t$ denotes its volume form. The factor $r(t)$ is used to normalize the equation so that the volume is constant.
 
 Assuming that the connection $A_0$ is Yang-Mills then, under the same setting as in Section \ref{sec21}, Corollary \ref{scalarwithfconst} implies that 
\[ R(t) = \pi^*R_{g_t}-\frac{1}{2}f_t^2\pi^*|F_{A_t}|^2_{g_t} = 2K_0\lambda_t^{-1}-\frac{1}{2}F_0f^2_t \lambda^{-2}_t. \]
We observe that the scalar curvature $R(t)$ is constant on $P$, thus $r(t)=R(t)$. Therefore, the normalized Ricci flow equation \eqref{NRFEq}) can be written as the system
\begin{align*}
\dot{\lambda_t}&= -\frac{2}{3}K_0+\frac{2}{3}F_0f^2_t\lambda_t^{-1},\quad \text{for } \lambda_0=1,\\ 
\dot{f_t}&= \frac{2}{3}K_0f_t\lambda_t^{-1}-\frac{2}{3}F_0f^3_t\lambda_t^{-2}, \quad \text{for } f_0=1.
\end{align*}

In the remaining part of the paper, we study the solutions of the above system considering all possible cases for the Gaussian curvature $K_0$ and the intensity $F_0$.
Using the fact that 
\[ \dot{\lambda_t}f_t+\lambda_t\dot{f_t}=0 \]
our system can be reduced to
\begin{equation}\label{RF1234}
\dot{\lambda_t}=-\frac{2}{3}K_0+\frac{2}{3}F_0\lambda_t^{-3} , \quad f_t=\lambda_t^{-1}, \quad \lambda_0 =f_0=1.
\end{equation}
Next, we consider the following cases depending on $K_0$ and $F_0$. 
\begin{enumerate}
\item When the connection $A_0$ is flat, that is, $F_0=0$ and $K_0 = 0$, \eqref{RF1234} gives 
\[ \lambda_t=f_t=1. \]
Hence, the Ricci flow solution is trivial, 
\[ g_A(t)=g_A(0), \]
and the circle bundle has a geometric structure with model geometry $\mathbf{E}^3$. 
\item When $F_0=0$ and $K_0\neq 0$, we obtain
\[ \dot{\lambda_t}=-\frac{2}{3}K_0 \] 
which gives
\[ \lambda_t=1-\frac{2}{3}K_0t. \]
Therefore, the solution of the Ricci flow equation is given by
\[ g_A(t)=\Big(1-\frac{2}{3}K_0t\Big)\pi^*g_0+\Big(1-\frac{2}{3}K_0t\Big)^{-2}\,\omega_0\otimes\omega_0 \]
and the scalar curvature of this metric is given by 
\[ R_{t}=6K_0(3-2K_0t)^{-1}. \]
Note that if $K_0 > 0$, a singularity appears both in the metric and in the curvature at time $t=3/(2K_0)$. In particular, the metric on the base vanishes and the fiber expands. In this case, the circle bundle has a geometric structure with model geometry $\mathbf{S}^2\times \mathbf{R} $. Interestingly, this is the only class of locally homogeneous geometries whose Ricci flow achieve curvature singularities (see \cite{IJ}).

On the other hand, if $K_0 <0$, then the Ricci flow solution exists for any time $t\geq0$. If $t$ goes to infinity, then $\lambda_t\to \infty$ and $f_t\to 0$, which means that the metric on the base blows up and the fiber collapses. Moreover, the scalar curvature falls off to zero as well. In this case, the circle bundle has a geometric structure with model geometry $\mathbf{H}^2\times \mathbf{R} $.

\item When $F_0>0 $ and $K_0=0$, then \eqref{RF1234} takes the form
\[ \dot{\lambda_t}=\frac{2}{3}F_0\lambda_t^{-3}, \] 
which gives
\[ \lambda_t=\Big(1+\frac{8}{3}F_0 t\Big)^{1/4}. \]
Hence, the Ricci flow solution is given by 
\[ g_A(t)=\Big(1+\frac{8}{3}F_0t\Big)^{1/4}\pi^*g_0+\Big(1+\frac{8}{3}F_0t\Big)^{-1/2} \omega_0\otimes\omega_0. \]
Note that, as $t$ goes to infinity, then $\lambda_t\to\infty$ and $f_t\to 0$, which means that the metric $g_t$ on the base manifold $M$ expands and the fiber collapses. The scalar curvature of $g_A(t)$ is given by 
\[ R_{t}=-\frac{1}{2}F_0 \Big( 1+\frac{8}{3}F_0 t \Big)^{-1/2}. \]
Hence, the scalar curvature falls off to zero as $t$ goes to infinity. In this case, the circle bundle $P$ over the surface $M$ has a geometric structure with model geometry $\mathbf{Nil}$ (the geometry of the Heisenberg group).
\item When $F_0>0$ and $K_0<0$, we have 
\[ \dot{\lambda_t} =-\frac{2}{3}K_0+\frac{2}{3}F_0f_t^2\lambda_t^{-1} = -\frac{2}{3}K_0+\frac{2}{3}F_0\lambda_t^{-3} > \frac{2}{3}F_0\lambda_t^{-3}. \]
As a result, we obtain 
\[ \lambda_t^4> \frac{2}{3}F_0t+C_0, \quad C_0\in \mathbb{R}. \]
Therefore, $\lambda_t\to\infty$ and $f_t\to0$ as $t$ goes to infinity. This means that the metric on the base expands while the fiber of the bundle shrinks to zero. Moreover, we can see that the scalar curvature falls off to zero as $t\to\infty$. In this case, the circle bundle has a geometric structure with model geometry of the universal cover of $\mathrm{SL}(2, \R)$.

\item Finally, we consider $F_0>0$ and $K_0>0$. Then, from \eqref {RF1234} we have
\[ \dot{\lambda_t}=-\frac{2}{3}K_0-\frac{2}{3}F_0\lambda_t^{-3}. \]
Solving this equation leads to the implicit solution
\[ c \ln\frac{\lambda^2_t-c\lambda_t+c^2}{(\lambda_t+c)^2}- 2c \sqrt 3\tan^{-1}\bigg(\frac{2\lambda_t-c}{c\sqrt 3 }\bigg)+6\lambda_t=-4K_0 t +C, \]
where $c = -\sqrt[3]{F_0/K_0}$ and $C$ is a constant.
\end{enumerate}

One can see that as $t$ goes to $\infty$, the right-hand side goes to $-\infty$. The left-hand side goes to $-\infty$ for $\lambda_t\to -c$ and $f_t\to -1/c$. The scalar curvature of the bundle metric approaches to $2K_0+F_0/(2c^2)$ as $t$ goes to infinity. In this case, the circle bundle has a spherical geometry $\mathbf{S^3}$.

\section{Conclusion}\label{sec4}

In this paper, we considered a class of metrics on the total space of a principal bundle
over a Riemannian manifold defined using a connection on the
bundle and a metric on the base manifold. 
Specifically, we studied the Ricci flow equation on circle bundles over surfaces equipped with locally homogeneous connections. Through a comprehensive derivation and analysis of the Ricci flow equation, we demonstrated that the solutions correspond to six of the eight model geometries outlined in Thurston’s geometrization conjecture. 

The novelty of this work lies in the introduction of an alternative approach, wherein we simplify the Ricci flow equation into a system of nonlinear ordinary differential equations. This enables us to derive explicit solutions and analyze their asymptotic behavior, thereby deepening our understanding of the geometric structures on these manifolds

This work not only builds upon existing research in Ricci flow but also paves the way for future studies, such as the exploration of $\mathrm{SU}(2)$-bundles with locally homogeneous connections.

\vspace{3mm} \noindent
 {\bf Acknowledgments.}
 The authors sincerely thank the reviewers for their valuable feedback and constructive suggestions, which have greatly contributed to the improvement of the manuscript.
  The authors supported by a Xiamen University Malaysia Research Fund (Grant No: XMUMRF/2023-C12/IMAT/0027).

\vspace{3mm}  \noindent
{\bf Statements and Declarations.} The authors declare no competing interests.

\bibliographystyle{amsplain}
\bibliography{Bazdar-Fotopoulos}
\end{document}